\documentclass[11pt,bezier]{article}
\usepackage{amsmath,graphicx,amssymb,amsfonts,color}

\font\bg=cmbx10 scaled\magstep1
\font\Bg=cmbx12 scaled\magstep3
\font\small=cmr8

\newtheorem{newlemma}{{\bf Lemma}}

\newenvironment{lema}{\begin{newlemma}{\hspace{-0.5
em}{\bf.}}}{\end{newlemma}}

\newtheorem{newteorem}{{\bf Theorem}}

\newenvironment{teorem}{\begin{newteorem}{\hspace{-0.5
em}{\bf.}}}{\end{newteorem}}

\newtheorem{newkorolari}{{\bf Corollary}}

\newenvironment{korolari}{\begin{newkorolari}{\hspace{-0.5
em}{\bf.}}}{\end{newkorolari}}

\newtheorem{newdefine}{{\bf Definition}}

\newtheorem{newquestion}{{\bf Question}}

\newtheorem{newkonjek}{{\bf Conjecture}}

\newtheorem{newexample}{{\bf Example}}

\begin{document}
\tolerance=10000
\baselineskip18truept
\newbox\thebox
\global\setbox\thebox=\vbox to 0.2truecm{\hsize
0.15truecm\noindent\hfill}
\def\boxit#1{\vbox{\hrule\hbox{\vrule\kern0pt
     \vbox{\kern0pt#1\kern0pt}\kern0pt\vrule}\hrule}}
\def\qed{\lower0.1cm\hbox{\noindent \boxit{\copy\thebox}}\bigskip}
\def\ss{\smallskip}
\def\ms{\medskip}
\def\bs{\bigskip}
\def\c{\centerline}
\def\nt{\noindent}
\def\ul{\underline}
\def\ol{\overline}
\def\lc{\lceil}
\def\rc{\rceil}
\def\lf{\lfloor}
\def\rf{\rfloor}
\def\ov{\over}
\def\t{\tau}
\def\th{\theta}
\def\k{\kappa}
\def\l{\lambda}
\def\L{\Lambda}
\def\g{\gamma}
\def\d{\delta}
\def\D{\Delta}
\def\e{\epsilon}
\def\lg{\langle}
\def\rg{\rangle}
\def\p{\prime}
\def\sg{\sigma}
\def\ch{\choose}

\newcommand{\ben}{\begin{enumerate}}
\newcommand{\een}{\end{enumerate}}
\newcommand{\bit}{\begin{itemize}}
\newcommand{\eit}{\end{itemize}}
\newcommand{\bea}{\begin{eqnarray*}}
\newcommand{\eea}{\end{eqnarray*}}
\newcommand{\bear}{\begin{eqnarray}}
\newcommand{\eear}{\end{eqnarray}}

\centerline{\Bg The Domination Polynomials of  Cubic}
 \vspace{.3cm}

\centerline {\Bg  Graphs of Order 10}
\bigskip

\baselineskip12truept \centerline{\bg Saieed Akbari$^{a}$, Saeid
Alikhani$^{b,d,}${}\footnote{\baselineskip12truept\it\small
Corresponding author. E-mail: alikhani@yazduni.ac.ir}, Yee-hock
Peng$^{c,d}$} \baselineskip20truept \centerline{\it
$^{a}$Department of Mathematical Sciences, Sharif University of
Technology} \vskip-9truept \centerline{\it 11365-9415 Tehran,
Iran} \vskip-5truept \centerline{\it $^{b}$Department of
Mathematics, Yazd University } \vskip-9truept \centerline{\it
89195-741, Yazd, Iran} \vskip-5truept \centerline{\it
$^{c}$Department of Mathematics, and } \vskip-8truept
\centerline{\it $^{d}$Institute for Mathematical Research,
University Putra Malaysia} \vskip-9truept \centerline{\it  43400
UPM Serdang, Malaysia} \vskip-0.2truecm \nt\rule{16cm}{0.1mm}

\nt{\bg ABSTRACT}
\medskip

\baselineskip14truept

\nt{\it Let $G$ be a simple graph of order $n$.
 The domination polynomial of  $G$  is the polynomial $D(G,x)=\sum_{i=\gamma(G)}^{n} d(G,i) x^{i}$,
 where $d(G,i)$ is the number of dominating sets  of $G$ of size $i$, and
$\gamma(G)$ is the domination number of $G$. In this paper we study the
 domination polynomials of cubic graphs of order $10$.
 As a consequence, we show that the Petersen graph  is
  determined uniquely by its domination polynomial.}

\ms

\nt{\bf  Mathematics Subject Classification:} {\small 05C60.}

\nt\rule{16cm}{0.1mm}

\baselineskip20truept

\section{Introduction}

\nt Let $G=(V,E)$ be a simple graph.  The {\it order} of $G$ is the number of vertices of $G$.
 For any vertex $v\in V$, the {\it open neighborhood} of $v$ is the
 set $N(v)=\{u \in V|uv\in E\}$ and the {\it closed neighborhood} of $v$
is the set $N[v]=N(v)\cup \{v\}$. For a set $S\subseteq V$, the open
neighborhood of $S$ is $N(S)=\bigcup_{v\in S} N(v)$ and the closed neighborhood of $S$
 is $N[S]=N(S)\cup S$.
A set $S\subseteq V$ is a {\it dominating set} if $N[S]=V$, or equivalently,
every vertex in $V\backslash S$ is adjacent to at least one vertex in $S$.
The {\it domination number} $\gamma(G)$ is the minimum cardinality of a dominating set in $G$.
 A dominating set with cardinality $\gamma(G)$ is called a {\it $\gamma$-set}.
The family
of all $\gamma$-sets of a graph $G$ is denoted by $\Gamma(G)$.
For a detailed treatment of these parameters, the reader is referred to~\cite{domination}.
The {\it $i$-subset}  of $V(G)$ is a subset of $V(G)$  of size $i$.
Let ${\cal D}(G,i)$ be the family of dominating sets of a graph $G$ with cardinality $i$ and
 let $d(G,i)=|{\cal D}(G,i)|$.
 The {\it domination polynomial} $D(G,x)$ of $G$ is defined as
$D(G,x)=\sum_{ i=\gamma(G)}^{|V(G)|} d(G,i) x^{i}$,
where $\gamma(G)$ is the domination number of $G$ (see\cite{saeid1}).

\nt We  denote the family of all dominating sets
 of $G$ with cardinality  $i$ and contain a vertex $v$ by ${\cal D}_v(G,i)$, and $d_v(G,i)=|{\cal D}_v(G,i)|$.
Two graphs $G$ and $H$ are said to be {\it dominating equivalence}, or simply ${\cal D}$-equivalent,
written $G\sim H$, if $D(G,x)=D(H,x)$.
It is evident that the relation $\sim$ of being ${\cal D}$-equivalence
 is an equivalence relation on the family ${\cal G}$ of graphs, and thus ${\cal G}$ is partitioned into equivalence classes,
called the {\it ${\cal D}$-equivalence classes}. Given $G\in {\cal G}$, let
\[
[G]=\{H\in {\cal G}:H\sim G\}.
\]
We call $[G]$ the equivalence class determined by $G$.
A graph $G$ is said to be {\it dominating unique}, or simply {\it ${\cal D}$-unique}, if $[G]=\{G\}$.

\nt The minimum  degree of   $G$  is denoted by
$\delta(G)$. A graph $G$ is called {\it $k$-regular} if all
vertices  have the same degree $k$.
 A {\it vertex-transitive graph} is a graph $G$ such that for every pair of vertices $v$ and $w$ of $G$,
there exists an automorphism $\theta$ such that $\theta(v)=w$.
 One of the famous graphs is the Petersen
graph which is a symmetric non-planar cubic graph.
 In the study of domination polynomials, it is interesting to investigate
 the dominating sets and domination polynomial of this graph. We denote the Petersen graph by $P$.

\ms

\nt In this paper,  we study the dominating sets and domination polynomials of
cubic graphs of order 10. As a consequence, we show that the Petersen graph
is determined uniquely by its domination polynomial.
In the next section, we obtain domination polynomial of the  Petersen graph.
In Section 3, we list all $\gamma$-sets of connected cubic graphs of order 10.
This list will be used to study the ${\cal D}$-equivalence of these graphs in the last section.
In Section 4, we prove that
 the Petersen graph is ${\cal D}$-unique. In the last section, we study
  the ${\cal D}$-equivalence classes
 of some cubic graphs of order $10$.

\section{Domination Polynomial of the Petersen Graph}

\nt In this section we shall investigate the domination polynomial of  the Petersen graph.
 First, we state and prove  the following lemma:

\begin{lema}\label{lemma1}
Let $G$ be a vertex transitive graph of order $n$ and $v\in
V(G)$. For any $ 1\leq i\leq n$,  $d(G, i)=\frac{n}{i}d_v(G,i)$.
\end{lema}
\nt{\bf Proof.}
Clearly, if $D$ is a dominating set of $G$ with size $i$, and $\theta\in
Aut(G)$, then $\theta(D)$ is also a dominating set of $G$ with size $i$.
Since $G$ is a vertex transitive graph, then for every two vertices $v$
and $w$, $d_v(G,i)=d_w(G,i)$.
If $D$ is a dominating set of size $i$, then there are exactly $i$
vertices  $v_{j_1}, \ldots, v_{j_i}$ such that $D$ counted in
$d_{v_{j_r}}(G, i)$, for each $1\leq r\leq i$. Hence $d(G,
i)=\frac{n}{i}d_v(G,i)$, and the proof is complete.\quad\qed

\begin{teorem}\label{theorem1} {\rm(\cite{domination}, p.48)}
If $G$ is a connected graph of order $n$ with $\delta(G)\geq 3$,
 then $\gamma(G)\leq\frac{3n}{8}$.
\end{teorem}

\nt We need the following theorem for finding the domination
polynomial of the Petersen graph.

\begin{teorem}\label{theorem2}
Let $G$ be a graph of order $n$ with domination polynomial
$D(G,x)=\sum_{i=1}^n d(G,i)x^i$. If $d(G,j)={n \choose j}$ for
some $j$, then $\delta(G)\geq n-j$.  More precisely,
$\delta(G)=n-l$, where $l=\min \Big\{j|d(G,j)={n \choose
j}\Big\}$, and there are at least ${n \choose l-1}-d(G,l-1)$
vertices of degree $\delta(G)$ in $G$. Furthermore, if for every
two vertices of degree $\delta(G)$, say $u$ and $v$ we have
$N[u]\neq N[v]$, then there are exactly ${n \choose
l-1}-d(G,l-1)$ vertices of degree $\delta(G)$.
\end{teorem}

\nt{\bf Proof.} Since $d(G,j)={n \choose j}$ for any  $r\geq j$, every $r$-subset
of the vertices of  $G$ forms a dominating set for $G$. Suppose that
there is a vertex $v\in V(G)$ such that $deg(v)< n-j$.
Consider $V(G)\backslash N[v]$.
Clearly, $\Big|V(G)\backslash N[v]\Big|\geq j$, and $V(G)\backslash N[v]$ is not a
dominating set for $G$, a contradiction. Now, assume that $S$ is a
$(l-1)$-subset of $V(G)$ which is not a dominating set. Thus there is a vertex
$u\in V(G)\backslash S$ which is not covered by $S$.
 Since $\delta(G)\geq n-l$, we have $deg(u)=n-l$.
Let $S'\neq S$ be a $(l-1)$-subset which is not a dominating set.
Thus there exists a vertex
$u'\in V(G) \backslash S'$ which is not covered by $S'$. As we did before
$deg(u')=n-l$. We claim that $u\neq u'$.
Since $S\neq S'$, there exists a vertex $x\in S\cap (V(G)\backslash S')$.
 We know that $u'x\in E(G)$.
If $u=u'$, then $ux\in E(G)$, a contradiction.
 Thus $u\neq u'$ and the
claim is proved. If $u$ and $v$ are two vertices of degree
$\delta(G)$ and $N[u]=N[v]$, then $V(G)\backslash
N[u]=V(G)\backslash N[v]$ is an $(l-1)$-subset which is not a
dominating set for $G$. Thus we have at least ${n \choose
l-1}-d(G,l-1)$  vertices of degree $\delta(G)$ in $G$.  The last
part of theorem is obvious. \quad\qed

\begin{figure}[h]
\hglue2.5cm
\includegraphics[width=11cm,height=5.1cm]{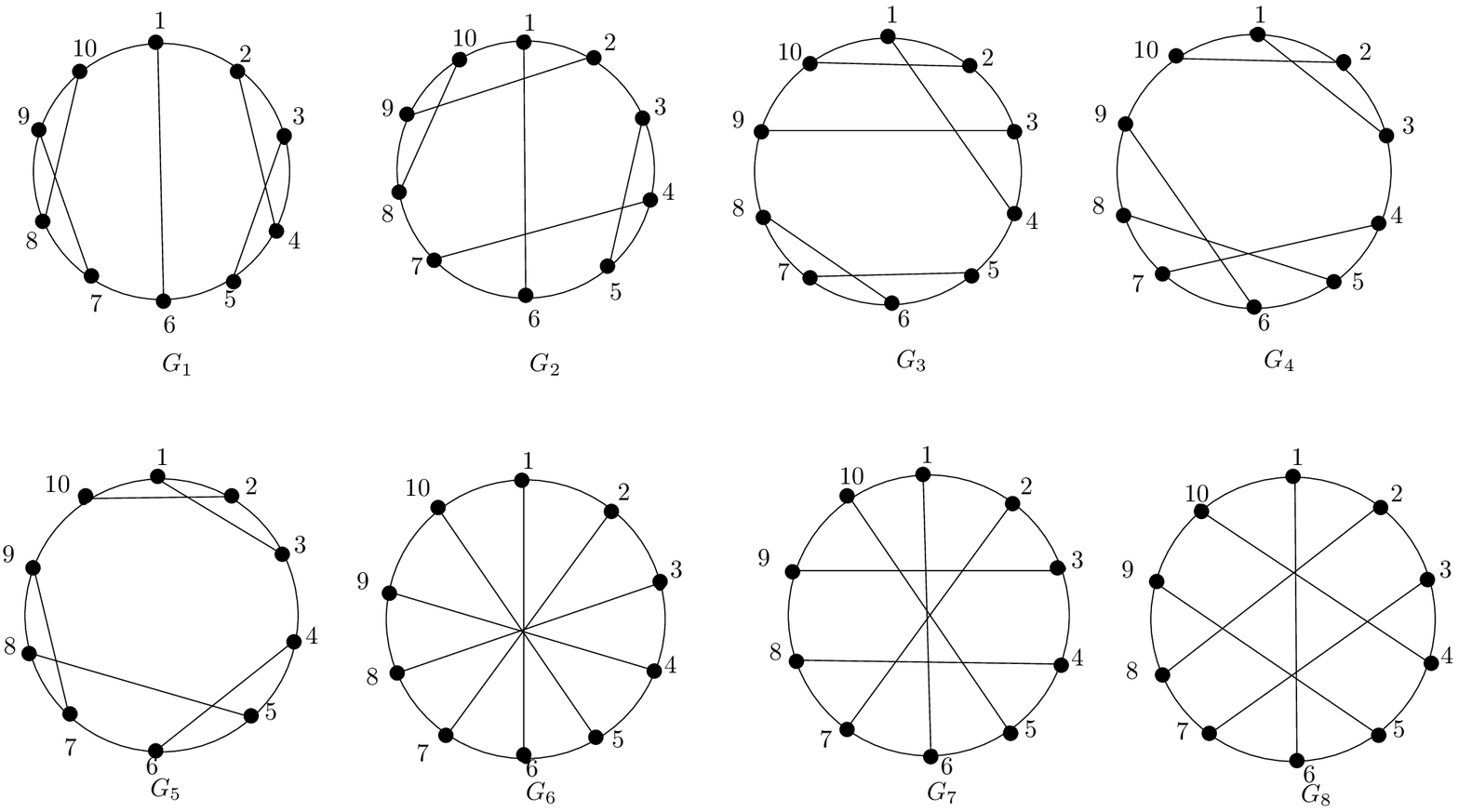}
\vglue5pt
\hglue2.5cm
\includegraphics[width=11cm,height=5.1cm]{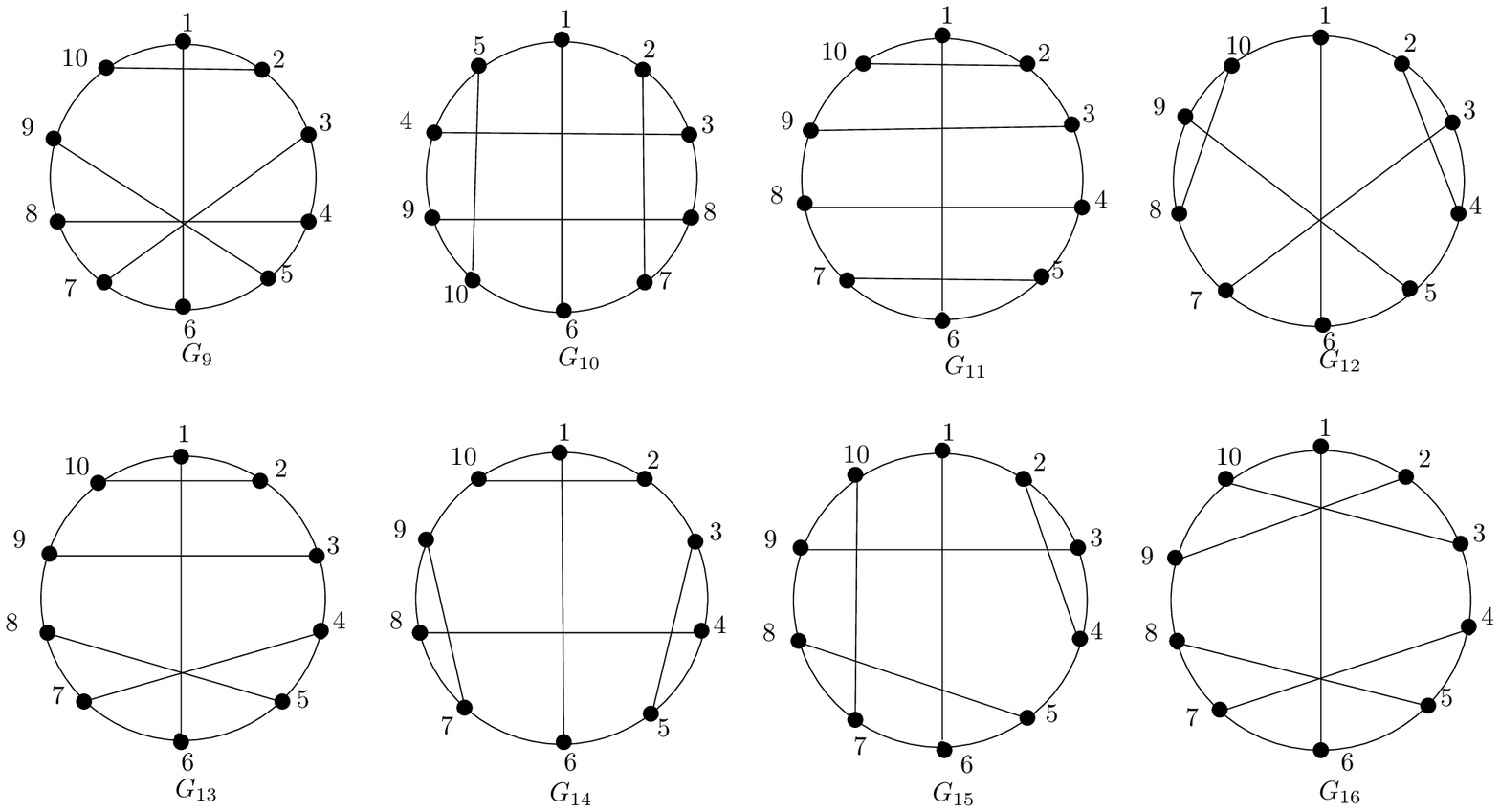}
\hglue2.5cm
\includegraphics[width=10.7cm,height=5cm]{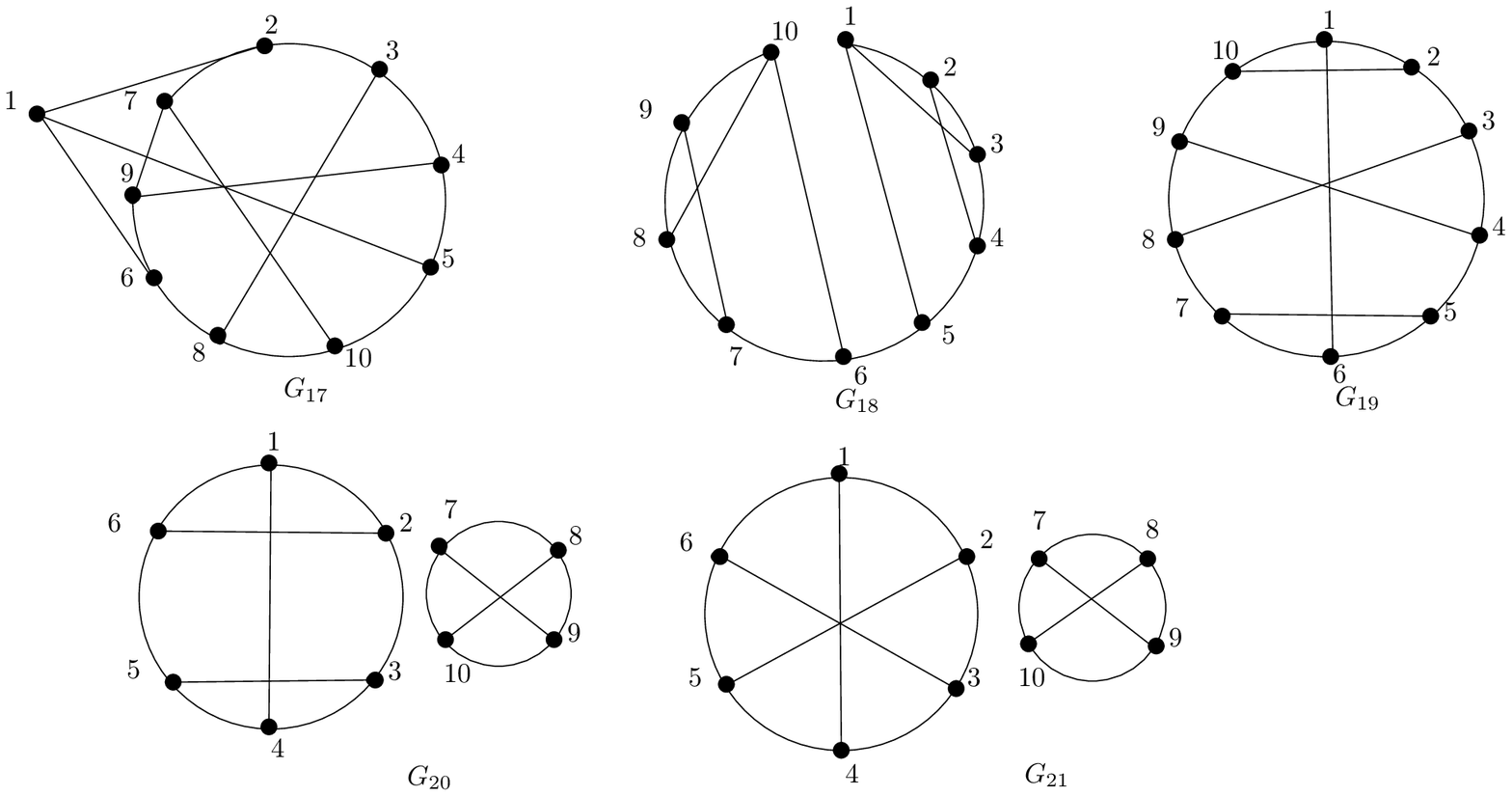}
\vglue-10pt \caption{\label{figure2} Cubic graphs of order 10.}
\end{figure}

\nt Indeed, by Theorem~\ref{theorem2},
 we have the following theorem which relates
the domination polynomial  and the regularity of a graph  $G$.

\begin{teorem}\label{theorem4}
Let $H$ be a $k$-regular graph, where for every two vertices $u,v\in
V(H)$, $N[u]\neq N[v]$. If  $D(G,x)=D(H,x)$,
then $G$ is also a $k$-regular graph.
\end{teorem}

\nt  There are exactly 21 cubic graphs of order $10$ given in
Figure~\ref{figure2} (see~\cite{reza}).
  Using Theorem~\ref{theorem1}, the domination number of a connected cubic
 graph of order $10$ is $3$. Three are just two non-connected cubic graphs of order $10$.
  Clearly, for  these graphs,  the domination number is also 3. Note that the graph $G_{17}$ is
 the Petersen graph.  For the labeled  graph $G_{17}$ in Figure~\ref{figure2}, we obtain  all dominating sets
  of size $3$ and $4$ in the following lemma.



\begin{lema} \label{lemma2}
For the Petersen graph $P$, $d(P,3)=10$ and $d(P,4)=75$.
\end{lema}
\nt{\bf Proof.} First, we list all dominating sets of $P$ of
cardinality $3$, which are the $\gamma$-sets of the labeled Petersen
graph (graph $G_{17}$)  given in Figure~\ref{figure2}.
\\
$ {\cal D}(P,3)=\Big\{\{1,3,7\},\{1,4,10\},\{1,8,9\},\{2,4,8\},
\{2,5,6\}, \{2,9,10\},
\{3,5,9\},\{3,6,10\},\{4,6,7\},\\\{5,7,8\}\Big\}.$ Now, we shall
compute $d(P,4)$. By Lemma~\ref{lemma1}, it suffices to obtain
the dominating sets of cardinality 4 containing one vertex, say
the vertex labeled $1$. These dominating sets are listed below.
\\
${\cal D}_1(P,4)=\Big\{\{1,2,3,7\}, \{1,2,4,8\},\{1,2,4,10\},\{1,2,5,6\},\{1,2,8,9\},
\{1,2,9,10\},\{1,3,4,7\},\\\{1,3,4,10\}, \{1,3,5,7\},\{1,3,5,9\},\{1,3,6,7\},\{1,3,6,10\},
\{1,3,7,8\},\{1,3,7,9\},\{1,3,7,10\},\\ \{1,3,8,9\},\{1,3,9,10\},\{1,4,5,10\},\{1,4,6,7\},
\{1,4,6,10\},\{1,4,7,8\}, \{1,4,7,10\}, \{1,4,8,9\},\\ \{1,4,8,10\}, \{1,4,9,10\},\{1,5,7,8\},
\{1,5,8,9\},\{1,6,8,9\},\{1,7,8,9\},\{1,8,9,10\}\Big\}$.

\nt Therefore by Lemma~\ref{lemma1}, $d(P,4)=\frac{10\times30}{4}=75$.\quad\qed

\begin{figure}[h]
\hspace{4cm}
\includegraphics[width=7.8cm,height=3.8cm]{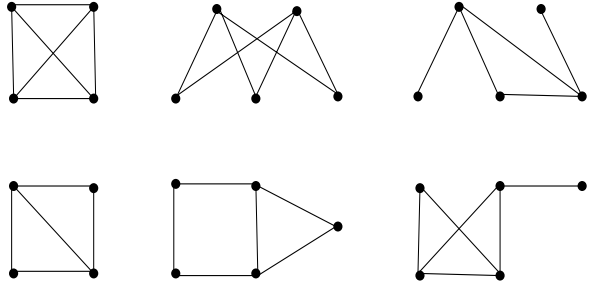}
\caption{\label{figure3} Graphs illustrated in Lemma~\ref{lemma5}.}
\end{figure}

\nt We need the following lemma:

\begin{lema}\label{lemma5}
Let $G$ be a cubic graph of order $10$. Then the following hold:
\begin{enumerate}

\item[(i)] $d(G, i)={n \choose i}$,  for $i=7,8,9,10$.

\item[(ii)]
if $t$ and $s$ are the number of subgraphs isomorphic to
$K_4\backslash\{e\}$($e$ is an edge) and $K_4$ in $G$,
respectively, then $d(G,6)={10 \choose 6}-(10-t-3s)$.

\item[(iii)]
if $G$ has no subgraph isomorphic to graphs given in Figure~\ref{figure3},
then $d(G,5)={10 \choose 5}-60$.
\end{enumerate}
\end{lema}

\nt{\bf Proof.} (i) It follows from Theorem~\ref{theorem2}.

 (ii)
 If $G$ is a cubic graph of order 10, then for
every $v\in V(G)$, $V(G)\backslash N[v]$ is not a dominating set.
Also, if $S\subset V(G)$, $|S|=6$ and $S$ is not a dominating set,
then $S=V(G)\backslash N[v]$, for some $v\in V(G)$. Note that if
$G$ has $K_4\backslash \{e\}$ as a subgraph, then there are two
vertices  $u_1$ and $u_2$ such that $G\backslash
N[u_1]=G\backslash N[u_2]$. Also if $G$ has $K_4$ as its subgraph,
then there are four vertices $u_i$, $1\leq i\leq 4$ such that
$G\backslash N[u_i]=G\backslash N[u_j]$, for $1\leq i\neq j\leq 4$.
Hence we have $d(G,6)={10 \choose 6}-(10-t-3s)$.

(iii) It suffices to determine the number of $5$-subsets which are
not dominating set. Suppose that $S\subseteq V(G)$, $|S|=5$, and
$S$ is not a dominating set for $G$. Thus there exists $v\in
V(G)$ such that $N[v]\cap S=\emptyset$.
 Now, note that for every $x\in V(G)$, $V(G)\backslash (N[x]\cup \{y\})$,
where $y\in V(G)\backslash N[x]$ is a 5-subset which is not a
dominating set for $G$. Also since none of the graphs given in
Figure~\ref{figure3} is a subgraph of $G$,
 for every two distinct vertices $x$ and $x'$
and any two arbitrary vertices $y\in V(G)\backslash N[x]$ and
$y'\in V(G)\backslash N[x']$, we have $V(G)\backslash (N[x]\cup
\{y\})\neq V(G)\backslash (N[x']\cup \{y'\})$. This implies that
the number of $5$-subsets of $V(G)$ which are not dominating sets
is $10\times 6=60$. So we have $d(G,5)={10 \choose
5}-60$.\quad\qed

\begin{korolari}\label{cor} For cubic graphs of order $10$, the following hold:

{\rm (i)} If $G\in \{G_{20}, G_{21}\},$ then $d(G, 6)={10 \choose
6}-7$.

{\rm (ii)} If $G\in \{G_1, G_{18}\}$, then $d(G, 6)={10 \choose
6}-8$.

{\rm (iii)} If $G\in \{G_3, G_{5}\}$, then $d(G, 6)={10 \choose
6}-9$.

{\rm (iv)} For each $i, 1\leq i\leq 21$, if $i\not \in \{1, 3, 5,
 18, 20, 21\}$, then $d(G_i, 6)={10 \choose
6}-10$.

{\rm (v)}  If $G\in \{G_6, G_7, G_8, G_{10}, G_{17}\},$ then
$d(G,5)={10 \choose 5}-60$.
\end{korolari}

\begin{teorem} \label{theorem3}
The domination polynomial of the Petersen graph $P$ is:
\[
D(P,x)=x^{10}+{10\choose 9}x^9+{10\choose 8}x^8+{10\choose 7}x^7+\Big[{10\choose 6}-10\Big]x^6+\Big[{10\choose 5}-60\Big]x^5+75x^4+10x^3.
\]
\end{teorem}
\nt{\bf Proof.} The result follows from Lemma~\ref{lemma2} and
Corollary~\ref{cor}.\quad\qed

\section{$\gamma$-Sets of cubic graphs of order $10$.}

\nt In this section, we present all $\gamma$-sets of connected cubic graphs
$G_1,G_2,\ldots,G_{18},G_{19}$ shown in
Figure~\ref{figure2}. The results here will be useful in studying the
${\cal D}$-equivalence of these graphs in the last section.

\nt $\Gamma(G_1)=\Big\{\{1,3,7\},\{1,3,8\},\{1,3,9\},\{1,4,7\},\{1,4,8\},
\{1,4,9\},\{1,5,7\},\{1,5,8\},\{1,5,9\},\\\{2,5,8\},\{2,5,9\},\{2,6,8\},
\{2,6,9\},\{2,6,10\},\{3,6,8\},\{3,6,9\},\{3,6,10\},\{3,7,10\},
\{4,6,8\},\\\{4,6,9\},\{4,6,10\},\{4,7,10\}\Big\}$. Therefore, $d(G_1,3)=\big|\Gamma(G_1)\big|=22$.


\ms
\nt $\Gamma(G_2)=\Big\{\{1,3,8\},\{1,4,8\},\{1,4,9\},\{1,4,10\},\{1,5,8\}, \{2,5,8\},
\{3,6,8\},\{3,6,9\},\{3,6,10\},\\\{3,7,10\},\{4,6,9\},\{5,6,9\}\Big\}.$
 Therefore $d(G_2,3)=\big|\Gamma(G_2)\big|=12$.

\ms

\nt $\Gamma(G_3)=\Big\{\{1,3,6\},\{1,3,7\},\{1,4,8\},\{1,5,9\},\{1,6,9\},\{1,7,9\}, \{2,3,6\},
\{2,3,7\},\{2,4,8\},\\\{2,5,8\},\{2,5,9\},\{3,6,10\},\{3,7,10\},\{4,6,10\},\{4,7,10\},\{4,8,10\},
\{5,9,10\}\Big\}.$

\nt Therefore $d(G_3,3)=\big|\Gamma(G_3)\big|=17$.

 \ms

 \nt \nt $\Gamma(G_4)=\Big\{\{1,5,6\},\{1,5,8\},\{1,6,7\},\{1,7,8\},\{2,4,9\},\{2,5,6\}, \{2,5,8\},
\{2,6,7\},\{2,7,8\},\\\{3,4,9\},\{3,6,9\},\{3,8,9\},\{4,5,10\},\{4,7,10\},\{4,9,10\}\Big\}.$
 Therefore $d(G_4,3)=\big|\Gamma(G_4)\big|=15$.

\ms

\nt $\Gamma(G_5)=\Big\{\{1,4,7\},\{1,4,8\},\{1,4,9\},\{1,5,7\},\{1,5,8\},\{1,5,9\}, \{1,6,7\},
\{1,6,8\},\{1,6,9\},\\\{2,4,7\},\{2,4,8\},\{2,4,9\},\{2,5,7\},\{2,5,8\},\{2,5,9\},\{2,6,7\},
\{2,6,8\},\{2,6,9\},\{3,4,9\},\\\{3,5,9\},\{3,6,9\},\{4,7,10\},\{4,8,10\},\{4,9,10\}\Big\}.$
 Therefore $d(G_5,3)=\big|\Gamma(G_5)\big|=24$.

\ms

\nt $\Gamma(G_6)=\Big\{\{1,4,7\},\{1,4,8\},\{1,5,8\},\{2,5,8\},\{2,5,9\},\{2,6,9\}, \{3,6,9\},
\{3,6,10\},\{3,7,10\},\\\{4,7,10\}\Big\}.$
 Therefore $d(G_6,3)=\big|\Gamma(G_6)\big|=10$.

\ms

\nt $\Gamma(G_7)=\Big\{\{1,4,8\},\{2,5,8\},\{2,5,9\},\{3,6,9\},\{3,7,10\},\{4,7,10\}\Big\}.$

\nt Therefore $d(G_7,3)=\big|\Gamma(G_7)\big|=6$

\ms

\nt
$\Gamma(G_8)=\Big\{\{1,3,9\},\{1,4,8\},\{1,5,7\},\{2,6,10\},\{3,6,9\},\{4,6,8\}\Big\}$.

\nt Therefore $d(G_8,3)=\big|\Gamma(G_8)\big|=6.$

\ms

\nt
$\Gamma(G_9)=\Big\{\{1,3,9\},\{1,4,8\},\{1,5,7\},\{2,5,7\},\{2,5,8\},\{2,6,8\},
\{3,6,9\},\{4,6,10\},\{4,7,10\},\\\{5,7,10\}\Big\}$.
 Therefore $d(G_9,3)=\big|\Gamma(G_9)\big|=10.$

\ms

\nt $\Gamma(G_{10})=\Big\{\{1,2,9\},\{1,5,8\},\{1,8,9\},\{2,3,10\},\{2,9,10\},\{3,4,6\},
\{3,6,10\},\{4,5,7\},\{4,5,7\},\\\{4,6,7\}\Big\}$.
 Therefore $d(G_{10},3)=\big|\Gamma(G_{10})\big|=10.$

\ms

\nt
$ \Gamma(G_{11})=\Big\{\{1,3,7\},\{1,5,9\},\{2,3,7\},\{2,5,8\},\{2,5,9\},\{2,6,8\},\{2,7,8\},
\{3,6,9\},\{3,7,10\},\\\{4,5,10\},\{4,6,10\},\{4,7,10\}\Big\}$.
 Therefore $d(G_{11},3)=\big|\Gamma(G_{11})\big|=12$.

\ms

\nt
$\Gamma(G_{12})=\Big\{\{1,3,9\},\{1,4,8\},\{1,5,7\},\{2,5,8\},\{2,6,8\},\{2,6,9\},\{2,6,10\},\{2,7,9\},
\{3,5,10\},\\\{3,6,8\},\{3,6,9\},\{4,6,8\},\{4,6,9\},\{4,6,10\},\{4,7,10\}\Big\}.$
Therefore $d(G_{12},3)=\big|\Gamma(G_{12})\big|=15$.

\ms

\nt
$\Gamma(G_{13})=\Big\{\{1,3,8\},\{1,4,8\},\{1,4,9\},\{2,5,8\},\{2,7,8\},\{3,6,9\},
\{4,5,10\},\{4,7,10\}\Big\}$.
Therefore $d(G_{13},3)=\big|\Gamma(G_{13})\big|=8$.

\ms

\nt
$\Gamma(G_{14})=\Big\{\{1,3,7\},\{1,3,8\},\{1,3,9\},\{1,4,7\},\{1,4,8\},\{1,4,9\},\{1,5,7\},\{1,5,8\},\{1,5,9\},
\\\{2,3,7\},\{2,4,7\},\{2,5,7\},\{2,5,8\},\{2,5,9\},\{2,6,8\},\{3,6,9\},\{3,7,10\},\{4,6,10\},
\\\{4,7,10\},\{5,7,10\},\{5,8,10\},\{5,9,10\}\Big\}$.
Therefore $d(G_{14},3)=\big|\Gamma(G_{14})\big|=22$.

\ms

\nt
$\Gamma(G_{15})=\Big\{\{1,2,8\},\{1,3,8\},\{1,4,8\},\{2,5,10\},\{2,6,9\},\{2,7,8\}
\{3,5,10\},\{3,6,7\},\{3,6,9\},\\\{4,5,10\},\{4,6,9\},\{4,7,10\}\Big\}$.
Therefore $\big|\Gamma(G_{15})\big|=12$.

\ms

\nt
$\Gamma(G_{16})=\Big\{ \{1,3,8\},\{1,4,8\},\{1,4,9\},\{3,6,8\},\{3,6,9\},\{4,6,9\}\Big\}$.
Therefore $d(G_{16},3)=6$.

\ms

\nt
$\Gamma(G_{17})=\Big\{\{1,3,7\},\{1,4,10\},\{1,8,9\},\{2,4,8\}, \{2,5,6\}, \{2,9,10\},
\{3,5,9\},\{3,6,10\},\{4,6,7\},\\\{5,7,8\}\Big\}.$
Therefore, $d(G_{17},3)=10$.

\ms

\nt
$\Gamma(G_{18})=\Big\{\{1,5,8\},\{1,5,9\},\{2,5,8\}, \{2,5,9\}, \{2,6,7\},
\{2,6,8\},\{2,6,9\},\{2,6,10\},\{3,5,8\},\\\{3,5,9\},\{3,6,7\},\{3,6,8\},\{3,6,9\},
\{3,6,10\},\{4,5,8\},\{4,5,9\}\Big\}.$
Therefore, $d(G_{18},3)=16$.

\ms

\nt
$\Gamma(G_{19})=\Big\{\{1,4,8\},\{1,5,8\},\{2,4,7\},\{2,5,8\},\{2,5,9\},\{2,6,9\},\{2,7,9\},\\\{3,5,10\},\{3,6,9\},
\{3,6,10\},\{3,7,10\},\{4,7,10\},\{5,8,10\}\Big\}$. Therefore $\big|\Gamma(G_{19})\big|=13$.

\section{${\cal D}$-Equivalence class of the Petersen Graph}

\nt In this section we  show that the Petersen graph is ${\cal D}$-unique.

\begin{teorem}\label{theorem5}
The Petersen graph $P$ is ${\cal D}$-unique.
\end{teorem}
\nt{\bf Proof.}
 Assume that $G$ is a graph such that
$D(G,x)=D(P,x)$. Since for every two vertices $x,y\in V(P)$,
$N[x]\neq N[y]$, by Theorem~\ref{theorem4},  $G$ is a $3$-regular graph
of order $10$.
Using the $\Big|\Gamma(G_i)\Big|$ for $i=1,\ldots,21$ in Section 3 we reject some graphs from $[P]$.
 Since $d(G_6,3)=d(G_9,3)=d(G_{10},3)=d(G_{17},3)=10$, we  compare the cardinality of
 the families of dominating sets of these four graphs  of size $4$.
 \\
 ${\cal D}_1(G_6,4)=\Big\{\{1,2,3,4\},\{1,2,3,10\},\{1,2,4,7\},\{1,2,4,8\},\{1,2,4,9\},\{1,2,5,8\},\{1,2,5,9\},
 \\\{1,2,6,9\},\{1,2,9,10\}, \{1,3,4,6\},\{1,3,4,7\}, \{1,3,4,8\},
 \{1,3,5,8\},\{1,3,6,8\}, \{1,3,6,9\}, \\\{1,3,6,10\},\{1,3,7,10\},\{1,3,8,10\},\{1,4,5,7\},
 \{1,4,5,8\},\{1,4,6,7\},\{1,4,6,8\},\{1,4,6,9\},\\\{1,4,7,8\},\{1,4,7,9\}, \{1,4,7,10\},
 \{1,4,8,9\},\{1,4,8,10\},\{1,5,6,8\},\{1,5,7,8\},\{1,5,8,9\},\\\{1,5,8,10\}, \{1,6,8,9\},\{1,8,9,10\}\Big\}.$

 \nt Therefore, by Lemma~\ref{lemma1},  $d(G_6,4)=\frac{34\times 10}{4}=85>d(P,4)=75$.

 \nt Now, we obtain the family of all dominating sets of $G_{10}$ of size $4$.
 \\
 ${\cal D}_1(G_{10},4)=\Big\{\{1,2,3,9\},\{1,2,3,10\},\{1,2,4,9\},\{1,2,5,8\},\{1,2,5,9\},\{1,2,6,9\},\\\{1,2,7,9\},
 \{1,2,8,9\},\{1,2,9,10\},\{1,3,4,6\},\{1,3,5,8\},\{1,3,6,8\},\{1,3,6,9\},\{1,3,6,10\},\\\{1,3,7,9\},
 \{1,3,7,10\},\{1,3,8,9\},\{1,3,8,10\},\{1,4,5,7\},\{1,4,5,8\},\{1,4,6,7\},\{1,4,6,8\},\\\{1,4,6,9\},
 \{1,4,7,9\},\{1,4,7,10\},\{1,4,8,9\},\{1,4,8,10\},\{1,5,6,8\},\{1,5,7,8\},\{1,5,8,9\},\\
 \{1,5,8,10\},\{1,6,8,9\},\{1,7,8,9\},
 \{1,8,9,10\}\Big\}$.

  \nt Therefore, by Lemma~\ref{lemma1}, $d(G_{10},4)=\frac{34\times10}{4}=85>d(P,4)=75$.

 \nt Now, for graph $G_9$ we have,
 ${\cal D}(G_9,4)=
 \\
 \Big\{\{1,2,3,9\},\{1,2,4,8\},\{1,2,5,7\},\{1,2,5,8\},\{1,2,6,8\},\{1,2,8,9\},\{1,3,4,5\},
 \{1,3,4,8\},\\\{1,3,4,9\},\{1,3,4,10\},\{1,3,5,7\},\{1,3,5,8\},\{1,3,5,9\},\{1,3,6,8\},
 \{1,3,6,9\},\{1,3,7,9\},\\\{1,3,8,9\},\{1,3,9,10\},\{1,4,5,6\},\{1,4,5,7\},\{1,4,5,8\},
 \{1,4,6,8\},\{1,4,6,9\},\{1,4,6,10\},\\\{1,4,7,8\},\{1,4,7,9\},\{1,4,7,10\},\{1,4,8,9\},
 \{1,4,8,10\},\{1,5,6,7\},\{1,5,7,8\},\{1,5,7,9\},\\\{1,5,7,10\},\{1,6,7,8\},\{1,7,8,9\}
 \{2,3,4,5\},\{2,3,5,7\},\{2,3,5,8\},\{2,3,5,9\},\{2,3,6,8\},\\\{2,3,6,9\},
 \{2,3,7,9\},\{2,4,5,6\},\{2,4,5,7\},\{2,4,5,8\},\{2,4,6,8\},\{2,4,6,9\},\{2,4,6,10\},\\\{2,4,7,8\},
 \{2,4,7,9\},\{2,4,7,10\},\{2,5,6,7\},\{2,5,6,8\},\{2,5,6,9\},\{2,5,7,8\},\{2,5,7,9\},
 \\\{2,5,7,10\},\{2,5,8,9\},\{2,5,8,10\},\{2,6,7,8\},\{2,6,8,9\},\{2,6,8,10\},\{2,7,8,9\},\{3,4,5,10\},
 \\\{3,4,6,10\},\{3,4,6,10\},\{3,4,7,10\},\{3,5,6,9\},\{3,5,7,10\},\{3,5,8,10\},\{3,5,9,10\},\{3,6,7,10\},
\\ \{3,6,8,10\},\{3,6,7,9\},\{3,6,8,9\},\{3,6,9,10\},\{3,7,9,10\},\{4,5,6,10\},\{4,5,7,10\},\{4,5,8,10\},\\\{4,6,7,10\},
 \{4,6,8,10\},\{4,6,9,10\},\{4,7,8,10\},\{4,7,9,10\},\{5,6,7,10\},\{5,7,8,10\},\{5,7,9,10\},\\\{6,7,8,10\},
 \{6,7,8,10\},\{7,8,9,10\}\Big\}$.

\nt Therefore $d(G_9,4)=91>d(P,4).$

\nt  Hence $[P]=\{P\}$, and so the Petersen graph is ${\cal D}$-unique.\quad\qed

\nt By the arguments in the proof of Theorem~\ref{theorem5}, we have the following corollary.

\begin{korolari}
{\rm (i)}
 The graph $G_9$ is ${\cal D}$-unique,

{\rm (ii)} $[G_6]=\Big\{G_6,G_{10}\Big\}$
  with the following domination polynomial:
\[
x^{10}+{10\choose 9}x^9+{10\choose 8}x^8+{10\choose 7}x^7+({10\choose 6}-10)x^6+({10\choose 5}-60)x^5+85x^4+10x^3.
\]
\end{korolari}

\section{${\cal D}$-equivalence class of cubic  graphs of order 10}

\nt In this section, we shall study the ${\cal D}$-equivalence
classes of other cubic graphs of order $10$.

\nt We need the following theorem:
\begin{teorem}\label{theorem6}{\rm (\cite{saeid1})}
If a graph $G$ has  $m$ components $G_1,\ldots,G_m$, then $D(G,x)=D(G_1,x)\cdots D(G_m,x)$.\quad\qed
\end{teorem}

\begin{korolari}
Two graphs $G_{20}$ and $G_{21}$ are ${\cal D}$-equivalence, with the following domination polynomial:
\[
D(G_{20},x)=D(G_{21},x)=x^{10}+10x^9+45x^8+120x^7+203x^6+216x^5+134x^4+36x^3.
\]
\end{korolari}
\nt{\bf Proof.} Two graphs $G_{20}$ and $G_{21}$ are disconnected with two components. In other words
$G_{20}=H\cup K_4$ and $G_{21}=H'\cup K_4$, where $H$ and $H'$ are graphs with 6 vertices. It is not hard to see that
\[
D(H,x)=D(H',x)=x^6+6x^5+15x^4+20x^3+9x^2.
\]
On the other hand, $D(K_4,x)=x^4+4x^3+6x^2+4x$. By Theorem~\ref{theorem6}, we have the result.\quad\qed

\begin{teorem}\label{theorem7}
The graphs $G_{12}$,$G_{13}$,$G_{14}$,$G_{16}$, and $G_{19}$ in Figure~\ref{figure2} are ${\cal D}$-unique.
\end{teorem}
\nt{\bf Proof.}
\nt Using $\gamma$-sets in Section 3,  $\big|\Gamma(G_{12})\big|=15, \big|\Gamma(G_{13})\big|=7,
\big|\Gamma(G_{14})\big|=22$, and $\big|\Gamma(G_{19})\big|=13$. By comparing these numbers with the cardinality of
$\gamma$-sets of other 3-regular graphs, we have the result. Now,  we consider graph $G_{16}$.
 Since $d(G_7,3)=d(G_8,3)=d(G_{16},3)=6$, we shall obtain $d(G_i,4)$ for $i=7,8,16$.
 \\
 ${\cal D}_1(G_7,4)=\Big\{\{1,2,3,4\},\{1,2,4,8\},\{1,2,4,9\},\{1,2,4,10\},\{1,2,5,8\},\{1,2,5,9\},
 \{1,2,6,8\},\\\{1,2,8,10\},\{1,3,4,6\},\{1,3,4,7\},\{1,3,4,8\},\{1,3,5,7\},\{1,3,5,8\},\{1,3,6,7\},
 \{1,3,6,8\},\\\{1,3,6,9\},\{1,3,7,10\},\{1,3,8,10\},\{1,4,5,8\},\{1,4,6,8\},\{1,4,6,9\},\{1,4,6,10\},
 \{1,4,7,8\},\\\{1,4,7,9\},\{1,4,7,10\},\{1,4,8,9\},\{1,4,8,10\},\{1,5,6,9\},\{1,5,7,9\},\{1,5,8,9\},
 \{1,6,8,9\},\\\{1,8,9,10\}\Big\}$. Therefore, by Lemma~\ref{lemma1},  $d(G_7,4)=\frac{32\times 5}{2}=80$.

\begin{figure}[h]
\hspace{1.4cm}
\includegraphics[width=13cm,height=12cm]{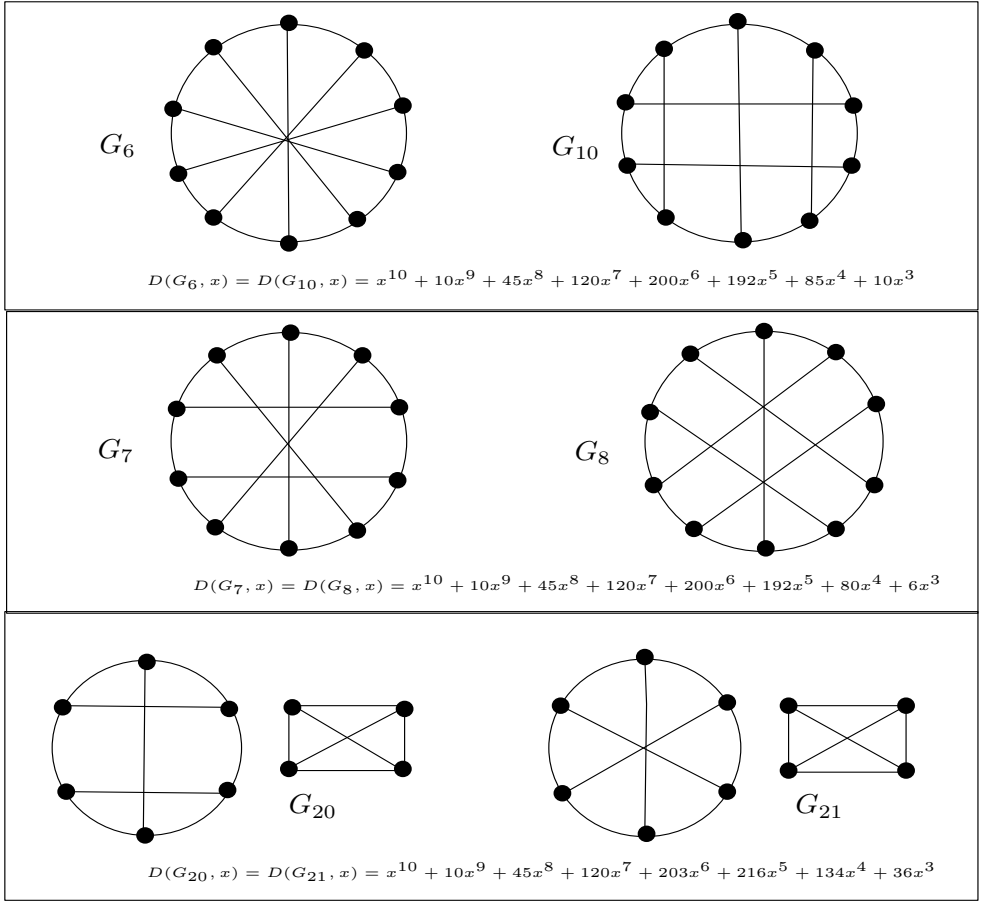}
\caption{\label{figure4} Cubic graphs of order 10 with identical domination polynomial.}
\end{figure}

\nt
 ${\cal D}_1(G_8,4)=\Big\{\{1,2,3,5\},\{1,2,3,9\},\{1,2,4,8\},\{1,2,5,6\},\{1,2,5,7\},\{1,2,5,8\},
 \{1,2,6,10\},\\\{1,3,4,8\},\{1,3,4,9\},\{1,3,5,7\},\{1,3,5,8\},\{1,3,5,9\},\{1,3,6,8\},\{1,3,6,9\},
 \{1,3,7,9\},\\\{1,3,8,9\},\{1,3,9,10\},\{1,4,5,7\},\{1,4,5,8\},\{1,4,6,8\},\{1,4,6,9\},\{1,4,7,8\},
 \{1,4,7,9\},\\\{1,4,7,10\},\{1,4,8,9\},\{1,4,8,10\},\{1,5,6,7\},\{1,5,7,8\},\{1,5,7,9\},\{1,5,7,10\},
 \{1,6,7,10\},\\\{1,7,9,10\}\Big\}$. Therefore, by Lemma~\ref{lemma1}, $d(G_8,4)=\frac{32\times5}{2}=80$.

 \nt Now, we obtain $d(G_{16},4)$.
 \\
 ${\cal D}_1(G_{16},4)=\Big\{\{1,2,3,8\},\{1,2,4,5\},\{1,2,4,7\},\{1,2,4,8\},\{1,2,4,9\},\{1,2,5,6\},
 \{1,2,5,7\},\\\{1,2,5,8\},\{1,2,6,7\},\{1,2,7,8\},\{1,3,4,8\},\{1,3,4,9\},\{1,3,5,8\},\{1,3,6,8\},
 \{1,3,6,9\},\\\{1,3,7,8\},\{1,3,8,9\},\{1,3,8,10\},\{1,4,5,8\},\{1,4,5,9\},\{1,4,5,10\},\{1,4,6,8\},
 \{1,4,6,9\},\\\{1,4,7,8\},\{1,4,7,9\},\{1,4,7,10\},\{1,4,8,9\},\{1,4,8,10\},\{1,4,9,10\},\{1,5,6,10\},
 \{1,5,7,10\},\\\{1,5,8,10\},\{1,6,7,10\},\{1,7,8,10\}\Big\}$.
 Therefore, by Lemma~\ref{lemma1}, $d(G_{16},4)=\frac{34\times 5}{2}=85$. Hence $[G_{16}]=\{G_{16}\}$.\quad\qed

\nt By the arguments in the proof of Theorem~\ref{theorem7}, we have the following corollary.

\begin{korolari}
Two graphs $G_{7}$ and $G_{8}$ are ${\cal D}$-equivalence.
\end{korolari}


\ms

\nt In summary, in this paper we showed that the Petersen graph is ${\cal D}$-unique. Also, we proved that
 the graphs $G_2,G_9,G_{11},G_{12},G_{13},
G_{14},G_{15},G_{16},G_{17}$, and $G_{19}$ are ${\cal D}$-unique, and
$[G_6]=\{G_6,G_{10}\}$, $[G_7]=\{G_7,G_8\}$, $[G_{20}]=\{G_{20},G_{21}\}$ (see Figure~\ref{figure4}).
We are not able to determine the ${\cal D}$-equivalence of $G_1$, $G_3$, $G_4$, $G_5$, and $G_{18}$,
but we think that they are ${\cal D}$-unique.

\nt {\bf Acknowledgement.}  The first author is indebted to the Institute for Mathematical
 Research (INSPEM) at University Putra Malaysia (UPM) for
the partial  support and hospitality during his visit.


\begin{thebibliography}{99}

\bibitem{saeid1}  S. Alikhani, Y. H. Peng, Introduction to Domination polynomial of a graphs, Ars Combinatoria, to appear.
\bibitem{domination} T.W. Haynes, S.T. Hedetniemi, P.J. Slater, Fundamentals of Domination in Graphs, Marcel Dekker, NewYork, 1998.

\bibitem{reza} G. B. Khosrovshahi, Ch. Maysoori, Tayfeh-Rezaie, A Note on 3-Factorizations of $K_{10}$, J. Combin. Designs 9 (2001), 379-383.
\end{thebibliography}
\end{document}